\documentclass[english]{smfart}

\usepackage[english]{babel}
\usepackage[T1]{fontenc}
\usepackage{cfr-lm}
\usepackage{microtype}
\usepackage{smfthm}
\usepackage{amssymb}
\usepackage{mathtools}
\usepackage[mathscr]{euscript}
\usepackage{enumitem}
\usepackage{stmaryrd}
\usepackage[dvipsnames]{xcolor}
\usepackage{hyperref}

\renewcommand{\bfseries}{\sbweight}
\DeclareMathAlphabet{\mathbf}{\encodingdefault}{\familydefault}{sb}{n}

\author{Téofil Adamski}
\email{teofil.adamski@univ-smb.fr}
\address{Université Savoie Mont Blanc, CNRS, LAMA, 73000 Chambéry, France}
\title{Non-Archimedean and motivic stationary phase formulas}

\hypersetup{
    colorlinks=true,
    linkcolor=Red,
    citecolor=Red,
    bookmarksdepth=3,
}

\usepackage[backend=biber, maxnames=50]{biblatex}
\usepackage{csquotes}
\renewbibmacro{in:}{}
\bibliography{main.bib}

\newcommand{\RR}{\mathbf{R}}
\newcommand{\CC}{\mathbf{C}}
\newcommand{\QQ}{\mathbf{Q}}
\newcommand{\ZZ}{\mathbf{Z}}
\newcommand{\FF}{\mathbf{F}}
\newcommand{\KK}{\mathbf{K}}
\newcommand{\LL}{\mathbf{L}}
\newcommand{\cf}{\hbox{\normalfont\textbf{\textl{1}}}}
\newcommand{\E}{\mathbf{E}}
\newcommand{\e}{\mathbf{e}}

\newcommand{\Field}{\mathbf{Field}}
\newcommand{\Def}{\mathbf{Def}}
\newcommand{\RDef}{\mathbf{RDef}}
\newcommand{\RDefexp}{\mathbf{RDef}^{\,\text{exp}}}
\newcommand{\K}{\mathrm{K_0}}
\newcommand{\Pres}{\mathscr{P}}

\renewcommand{\L}{\mathscr{L}}
\newcommand{\B}{\mathrm{B}}
\renewcommand{\O}{\mathscr{O}}
\newcommand{\m}{\mathfrak{m}}
\newcommand{\M}{\mathscr{M}}
\renewcommand{\S}{\mathscr{S}}

\newcommand{\Cons}{\mathscr{C}}
\newcommand{\cons}{\mathrm{C}}
\newcommand{\expe}{^\text{exp}}
\newcommand{\Sch}{\mathscr{S}}
\newcommand{\Inte}{\mathscr{I}}
\newcommand{\inte}{\mathrm{I}}
\DeclareMathOperator{\supp}{supp}

\DeclareMathOperator{\ord}{ord}
\DeclareMathOperator{\ac}{\overline{ac}}
\DeclareMathOperator{\grad}{grad}
\DeclareMathOperator{\Hess}{Hess}
\DeclareMathOperator{\Jac}{Jac}
\DeclareMathOperator{\Crit}{Crit}

\newcommand{\fnc}[4]{\left\lvert\begin{aligned}#1&\longrightarrow#2\\#3&\longmapsto#4\end{aligned}\right.\kern-\nulldelimiterspace}
\newcommand{\fonc}[5]{#1\colon\fnc{#2}{#3}{#4}{#5}}
\newcommand{\dd}{\mathop{}\mathopen{}\mathrm{d}}
\renewcommand{\to}{\longrightarrow}

\newcommand{\abs}[1]{\lvert#1\rvert}
\newcommand{\inner}[2]{\langle#1, #2\rangle}
\newcommand{\ls}[1]{(\!(#1)\!)}
\newcommand{\fs}[1]{\llbracket#1\rrbracket}

\allowdisplaybreaks

\makeatletter
\DeclareRobustCommand{\@pointrait}{.}
\SetSymbolFont{stmry}{bold}{U}{stmry}{m}{n}
\g@addto@macro\bfseries{\boldmath}

\def\th@plain{%
    \let\thm@indent\noindent
    \thm@headfont{\bfseries\smf@boldmath\itshape}%
    \thm@notefont{\bfseries\smf@boldmath\upshape}%
    \thm@preskip.5\linespacing \@plus .5\linespacing
    \thm@postskip\thm@preskip
    \thm@headpunct{\MakePointrait}
    \itshape}

\begin{document}

\frontmatter

\begin{abstract}
    In this article, for a non degenerate singular phase, we reconsider a stationary phase formula of Heifetz~\cite{Heifetz} in the non-Archimedean local field setting and give a motivic analogue using Cluckers--Loeser's motivic integration.
\end{abstract}

\maketitle

\tableofcontents

\mainmatter

\section*{Introduction}

In real analysis, the stationary phase method gives an asymptotic development of integrals of the type
\begin{equation}\label{eq:osc-int}
    \int_{\RR^n} \varphi(x)\exp(i\lambda f(x)) \dd x
\end{equation}
when the real parameter $\lambda$ goes to infinity and $i^2 = -1$. These integrals are called \emph{oscillatory integrals}. This method goes back to Laplace and well-known formulas are given in~\cite[§7.7]{Hormander}, for instance when the function~$\varphi$ is compactly supported and the function $f$ is smooth or with nondegenerate critical points. In the latter case, the asymptotic development of the integral~\eqref{eq:osc-int} is given by considering the function $f$ on a neighborhood of the critical locus and its complement. The proof of the stationary phase method can be deduced from the Riemann--Lebesgue lemma.

For instance, this method allows to compute a whole family of integrals, like Fourier transforms (i.e. the case where $f(x) = \inner x\xi$ for some vector $\xi$) or Bessel integrals. Moreover, it can be used in microlocal analysis (and more generally for partial differential equations) to study wave front sets (e.g. see~\cite[Theorem~8.1.9]{Hormander}) or to construct pullbacks of distributions (see~\cite[Theorem~8.2.4]{Hormander}). We can also refer to Sawyer in~\cite{Sawyer} for some applications with random walks and to Huxley in~\cite{Huxley} for some asymptotic expansions of exponential sum in analytic number theory.

This kind of integrals has an analogue in the $p$-adic setting, using the Haar measure on the locally compact group $\QQ_p$ and a continuous group morphism from~$\QQ_p$ to $\CC^\times$ to mimic the exponential function. In 1985, Heifetz established in~\cite{Heifetz} the $p$-adic version of the mentioned stationary phase formulas in the real context. Also, Igusa in~\cite[§10.2]{Igusa} or~\cite[§10]{Igusa1994} has developed another stationary phase formula which computes~$p$-adic local zeta functions (see~\cite{Veys} for a survey).

\bigbreak

We start, in §\ref{sec:nonarchimedean} (see Theorem~\ref{thm:FPS-p}), by reconsidering that $p$-adic version of Heifetz when the phase $f$ has nondegenerate critical points (see~\cite[Proposition~1.2]{Heifetz}) and for a non-Archimedean local field (e.g. fields like $\QQ_p$ or $\FF_p\ls t$), giving a more detailed proof and using the non-Archimedean version of the Morse lemma by Cluckers--Herremans in~\cite{CluckersHerremans} (see also~\cite[Proposition~2.5]{CluckersKollarMustata}).

In §\ref{sec:motivic}, we consider the case of the non-Archimedean field $k\ls t$ for a characteristic zero field $k$. In that setting, we need to use motivic integration introduced by Kontsevich, developed successively by Denef--Loeser~\cite{DenefLoeser1999}, Cluckers--Loeser~\cite{CluckersLoeser2008, CluckersLoeser2010} and Hrushovski--Kazhdan~\cite{HrushovskiKazhdan}. In particular, using parameters integrals with exponential of Cluckers--Loeser's motivic integration and a motivic version of the non-Archimedean Morse lemma, we set in Theorem~\ref{thm:FPS-mot} a motivic stationary phase formula which is the analogue of Theorem~\ref{thm:FPS-p}. This formula could be used in the study of motivic wave front sets of distributions~\cite{Raibaut} as in the real setting.

\section{The non-Archimedean local field case}
\label{sec:nonarchimedean}

\subsection{Integration and Fourier transform}

\subsubsection{Setting and integration}

We consider a non-Archimedean local field $\KK$, i.e. either
\begin{itemize}
    \item a finite field extension of the $p$-adic numbers field $\QQ_p$ for some prime number~$p$,

    \item or a field isomorphic to the field $K\ls t$ of Laurent series with coefficients in some finite field $K$.
\end{itemize}
We denote by $\ord \colon \KK \to \ZZ \cup \{+\infty\}$ its valuation, $\O_\KK$ its valuation ring, $\m_\KK$ the maximal ideal of the ring $\O_\KK$ and $k_\KK \coloneq \O_\KK/\m_\KK$ the residual field of the valued field~$\KK$ with $q_\KK$ elements and characteristic $p_\KK$. The ideal $\m_\KK$ is generated by an element~$\varpi_\KK$ of $\O_\KK$. We consider the absolute value $\abs{\cdot}$ on~$\KK$ defined by the equality
\begin{equation}\label{eq:absolute-value}
    \abs x \coloneq q_\KK^{-\ord x} \in \RR
\end{equation}
for all elements $x$ of $\KK$. The abelian group $(\KK, +)$ is locally compact for the valuation topology and we consider its Haar measure $\dd x$ such that
\[ \int_{\O_\KK} \dd x = 1. \]

\subsubsection{Fourier transform}
\label{sssec:Fourier}

Let $\psi \colon \KK \to \CC^\times$ be an additive character ---~namely a continuous group morphism~--- which is nontrivial on $\O_\KK$ and trivial on $\m_\KK$. Let~$n$ be a nonnegative integer. For an integrable function $\varphi \colon \KK^n \to \CC$, its \emph{Fourier transform} is the function~$\hat\varphi \colon \KK^n \to \CC$ defined by the equality
\[ \hat\varphi(\xi) \coloneq \int_{\KK^n} \varphi(x)\psi(\inner x\xi) \dd x \]
for all elements $\xi$ of $\KK^n$ where the symbol $\inner\cdot\cdot$ denotes the canonical inner product on the vector space $\KK^n$.

Let $\Omega$ be an open set of $\KK^n$. An integrable function $\varphi \colon \Omega \to \CC$ is a \emph{Schwartz--Bruhat function} on~$\Omega$ if it is locally constant and its support
\[ \supp\varphi \coloneq \overline{\{x \in \Omega \mid \varphi(x) \ne 0\}} \subset \KK^n \]
is bounded. We denote by~$\Sch(\Omega)$ the set of Schwartz--Bruhat functions on $\Omega$.

\begin{theo}\label{thm:Fourier-inv}
    The Fourier transform $\varphi \longmapsto \hat\varphi$ induces a $\CC$-linear automorphism on the space $\Sch(\KK^n)$. More precisely, for all Schwartz--Bruhat functions $\varphi$ on~$\KK^n$ and for all elements $x$ of $\KK^n$, we have the equality
    \[ \hat{\hat\varphi}(x) = \varphi(-x). \]
\end{theo}

From Theorem~\ref{thm:Fourier-inv} and Fubini's theorem, we deduce in a straightforward way a Fourier--Plancherel type formula.

\begin{coro}\label{coro:Plancherel}
    Let $f$ and $g$ be two Schwartz--Bruhat functions on $\KK^n$. Then
    \[ \int_{\KK^n} \hat f(x)g(x) \dd x = \int_{\KK^n} f(x)\hat g(x) \dd x. \]
\end{coro}

\subsection{A Morse lemma}

Until the end, we suppose that the residual characteristic of the field $\KK$ is not equal to~$2$, i.e. $p_\KK \ne 2$. In this case, Cluckers--Herremans in~\cite{CluckersHerremans} prove a Morse lemma in the $p$-adic setting. We extend to the field $\KK$ some of their definitions. Let $n$ and $m$ be two positive integers.

\begin{defi}\label{def:analytic}
    Let $A$ be an open set of $\O_\KK^n$ and $B$ be an open set of $\O_\KK^m$.
    \begin{itemize}
        \item A function $A \to \O_\KK$ is \emph{globally analytic} if it is given by a formal power series with coefficients in $\O_\KK$, i.e. an element of $\O_\KK\fs{x_1, \dots, x_n}$, which converges on the open set $A$.

        \item A function $A \to \O_\KK^m$ is globally analytic if its components are globally analytic.

        \item A function $A \to B$ is a \emph{bianalytic function} if it is a globally analytic bijection whose inverse is globally analytic.
    \end{itemize}
    Let $f \colon A \to \O_\KK$ be a globally analytic function and $\alpha$ be a positive integer.
    \begin{itemize}
        \item The \emph{gradient} at some point $x$ of $A$ of the function $f$ is the $m$-uplet $\grad f(x)$ given by the equality
        \[ \grad f(x) \coloneq \left(\frac{\partial f}{\partial x_1}(x), \dots, \frac{\partial f}{\partial x_n}(x)\right) \]
        and its \emph{Hessian} at some point $x$ of $A$ is the matrix $\Hess f(x)$ of dimension~$n$ given by the equality
        \[ \Hess f(x) \coloneq \left(\frac{\partial^2f}{\partial x_i\partial x_j}(x)\right)_{1 \leqslant i, j \leqslant n}. \]

        \item A \emph{critical point} of the function $f$ is a point $x$ of $A$ such that
        \[ \grad f(x) = 0. \]

        \item A critical point $x$ of the function $f$ is \emph{nondegenerate} if
        \[ \det(\Hess f(x)) \ne 0. \]

        \item A \emph{critical point modulo $\varpi_\KK^\alpha$} of the function $f$ is a point $x$ of $A$ such that
        \[ \grad f(x) \equiv 0 \mod{\varpi_\KK^\alpha}. \]

        \item A critical point $x$ of the function $f$ is \emph{nondegenerate modulo $\varpi_\KK$} if
        \[ \det(\Hess f(x)) \not\equiv 0 \mod{\varpi_\KK}. \]
    \end{itemize}
\end{defi}

\begin{rema}\label{rem:gen-Kt}
    Definition~\ref{def:analytic} is also valid for valued fields of the form $\KK = K\ls t$ for a characteristic zero field $K$. In this case, we choose the natural valuation and the uniformizer $\varpi_{K\ls t} = t$. See Definitions~\ref{def:analytic-mot1} and~\ref{def:analytic-mot2} for a version which is uniform in the field $K$.
\end{rema}

Similarly to the $p$-adic Morse lemma~\cite[Lemma 3.3]{CluckersHerremans} (see also~\cite[Proposition~2.5]{CluckersKollarMustata}) with $p \ne 2$, we have the following result.

\begin{theo}\label{thm:Morse}
    Let $\alpha$ be a positive integer. Let $f \colon (\varpi_\KK^\alpha\O_\KK^{})^n \to \O_\KK^{}$ be a globally analytic function such that the origin is a critical point modulo $\varpi_\KK^\alpha$ which is nondegenerate modulo $\varpi_\KK$. Then there is a unique critical point $c$ in $(\varpi_\KK^\alpha\O_\KK^{})^n$ of the function~$f$. Moreover, there is a bianalytic function
    \[ T = (T_1, \dots, T_n) \colon (\varpi_\KK^\alpha\O_\KK^{})^n \to (\varpi_\KK^\alpha\O_\KK^{})^n \]
    with $T(0) = 0$ and elements $a_1$, \dots, $a_n$ of $\O_\KK^\times$ such that
    \[ f(x) = f(c) + \sum_{i = 1}^n a_iT_i(x)^2 \]
    for all elements $x$ of $(\varpi_\KK^\alpha\O_\KK^{})^n$.
\end{theo}

\begin{proof}
    One easily adapts the proof of~\cite[Lemma 3.3]{CluckersHerremans} by replacing the polydisc~$(p\ZZ_p)^n$ by $(\varpi_\KK^\alpha\O_\KK^{})^n$ and the set $p\ZZ_p$ by $\O_\KK$ and by using Lemma~\ref{lemma:CluckersHerremans2} instead of Lemma~\ref{lemma:CluckersHerremans} which are corollaries of the Igusa's non-Archimedean local inverse function theorem (see~\cite[§2.2]{Igusa} for more details).
\end{proof}

\begin{lemm}[Cluckers--Herremans {\cite[Lemma 3.2]{CluckersHerremans}}]\label{lemma:CluckersHerremans}
    Let $g \colon \O_\KK^n \to \O_\KK^n$ be a globally analytic function such that
    \[ \det(\Jac g(0)) \not\equiv 0 \mod{\varpi_\KK} \qquad \text{and} \qquad
       g(0) \equiv 0 \mod{\varpi_\KK} \]
    where the notation $\Jac g(0)$ refers to the jacobian matrix at the origin of the function~$g$. Then the restriction $g \colon (\varpi_\KK\O_\KK)^n \to (\varpi_\KK\O_\KK)^n$ is a bianalytic function.
\end{lemm}

From this lemma, we deduce immediately the following result.

\begin{lemm}\label{lemma:CluckersHerremans2}
    Let $\beta$ be a nonnegative integer. Let $g \colon (\varpi_\KK^\beta\O_\KK^{})^n \to (\varpi_\KK^\beta\O_\KK^{})^n$ be a globally analytic function such that
    \[ \det(\Jac g(0)) \not\equiv 0 \mod{\varpi_\KK} \qquad \text{and} \qquad
       g(0) \equiv 0 \mod{\varpi_\KK^{\beta + 1}}. \]
    Then the restriction $g \colon (\varpi_\KK^{\beta + 1}\O_\KK^{})^n \to (\varpi_\KK^{\beta + 1}\O_\KK^{})^n$ is a bianalytic function.
\end{lemm}

\begin{proof}
    The function
    \[ \fonc{\tilde g}{\O_\KK^n}{\O_\KK^n,}{x}{\varpi_\KK^{-\beta}g(\varpi_\KK^\beta x)} \]
    satisfies the condition of Lemma~\ref{lemma:CluckersHerremans}. Applying this lemma, it induces a bianalytic function $\tilde g \colon (\varpi_\KK\O_\KK)^n \to (\varpi_\KK\O_\KK)^n$. This concludes the lemma.
\end{proof}

We note that a critical point is a critical point modulo $\varpi_\KK^\alpha$ for all positive integers~$\alpha$. Then by translation, the following result is a consequence of Theorem~\ref{thm:Morse}.

\begin{coro}\label{coro:Morse}
    Let $\Omega$ be an open set of $\O_\KK^n$. Let $f \colon \Omega \to \O_\KK$ be a globally analytic function with a critical point $x_0$ of $\Omega$ which is nondegenerate modulo~$\varpi_\KK$. Let $\alpha$ be a nonnegative integer such that $B \coloneq x_0 + \varpi_\KK^\alpha\O_\KK^n \subset \Omega$. Then the point~$x_0$ is the only critical point of the function $f$ on the polydisc $B$. Moreover, there exist a bianalytic function
    \[ T = (T_1, \dots, T_n) \colon (\varpi_\KK^\alpha\O_\KK^{})^n \to (\varpi_\KK^\alpha\O_\KK^{})^n \]
    with $T(0) = 0$ and elements $a_1$, \dots, $a_n$ of $\O_\KK^\times$ such that
    \[ f(x) = f(x_0) + \sum_{i = 1}^n a_iT_i(x - x_0)^2 \]
    for all elements $x$ of $B$.
\end{coro}

\subsection{The non-Archimedean stationary phase formula}

Following notation of~§\ref{sssec:Fourier}, we prove the following lemmas which will be used in the proof of Theorem~\ref{thm:FPS-p}.

\begin{lemm}\label{lemma:SB-Psi}
    Let $\lambda$ be an element of $\KK$. Then the function
    \[ \fnc{\O_\KK}{\CC,}{u}{\psi(\lambda u^2)} \]
    is locally constant.
\end{lemm}

\begin{proof}
    If $\lambda = 0$, then the result is immediate. We suppose that $\lambda \ne 0$. We set
    \[ \beta \coloneq \max\mathopen{}\left\{0, \frac{1 - \ord\lambda}{\varpi_\KK}, \frac{1 - \ord\lambda}{2\varpi_\KK}\right\}\kern-\nulldelimiterspace. \]
    We verify that, for all elements $x$ of $\O_\KK$, the function under consideration in the lemma is constant on the ball $x + \varpi_\KK^\beta\O_\KK^{}$. For that, we use the fact that $\psi = 1$ on~$\varpi_\KK\O_\KK$.
\end{proof}

\begin{lemm}\label{lemma:Gauss}
    Let $\alpha$ be an integer. Let $a$ and $b$ be two elements of $\KK^\times$ such that
    \[ \ord b - \ord a \geqslant \alpha \qquad \text{and} \qquad
       2\ord b - \ord a \geqslant 1. \]
    Then
    \[ \int_{\varpi_\KK^\alpha\O_\KK^{}} \psi(ax^2 + bx) \dd x = \int_{\varpi_\KK^\alpha\O_\KK^{}} \psi(ax^2) \dd x. \]
\end{lemm}

\begin{proof}
    We have
    \begin{align*}
        \int_{\varpi_\KK^\alpha\O_\KK^{}} \psi(ax^2 + bx) \dd x &= \int_{\varpi_\KK^\alpha\O_\KK^{}} \psi\mathopen{}\left(a\left[x + \frac{b}{2a}\right]^2 - \frac{b^2}{4a}\right) \dd x \\
        &= \psi\mathopen{}\left(-\frac{b^2}{2a}\right) \int_{\varpi_\KK^\alpha\O_\KK^{}} \psi\mathopen{}\left(a\left[x + \frac{b}{4a}\right]^2\right) \dd x.
    \end{align*}
    As $p_\KK \ne 2$, we can write $\ord(b/2a) = \ord b - \ord a \geqslant \alpha$ and we get the bijection
    \[ \fnc{\varpi_\KK^\alpha\O_\KK^{}}{\varpi_\KK^\alpha\O_\KK^{},}{x}{x + b/2a} \]
    and, with this change of variables, we obtain
    \[ \int_{\varpi_\KK^\alpha\O_\KK^{}} \psi(ax^2 + bx) \dd x = \psi\mathopen{}\left(-\frac{b^2}{4a}\right) \int_{\varpi_\KK^\alpha\O_\KK^{}} \psi(ay^2) \dd y. \]
    Moreover, as $\ord(-b^2/4a) = 2\ord b - \ord a \geqslant 1$, we get $\psi(-b^2/4a) = 1$ and so the lemma.
\end{proof}

\begin{defi}\label{def:osc-int-p}
    Let $\Omega$ be an open set of $\O_\KK^n$. Let $\varphi$ be a Schwartz--Bruhat function on~$\Omega$ and $f \colon \Omega \to \O_\KK$ be a globally analytic function. For any element~$\lambda$ of $\KK$, we consider the integral
    \[ I_{f, \varphi}(\lambda) \coloneq \int_\Omega \varphi(x)\psi(\lambda f(x)) \dd x. \]
\end{defi}

We reconsider now~\cite[Proposition 1.2]{Heifetz} for a non-Archimedean local field. We add assumptions which seem to have been omitted in the statement of Heifetz and we give a detailed proof whose strategy will be used for the motivic case (see Theorem~\ref{thm:FPS-mot}).

\begin{theo}\label{thm:FPS-p}
    We use settings of Definition~\ref{def:osc-int-p} and make the following assumptions:
    \begin{enumerate}[label=(\roman*)]
        \item\label{it:xO} the function $f$ admits a unique critical point $x_0$ in $\Omega$;

        \item\label{it:xO-non-degen} the latter belongs to the set $\supp\varphi$ and it is nondegenerate modulo~$\varpi_\KK$;

        \item\label{it:def-R} there exists a function $R \colon \Omega \times \KK^n \to \KK$ such that
        \[ f(x + y) = f(x) + \inner{\grad f(x)}{y} + \inner{R(x, y)y}{y} \]
        for all elements $x$ of $\Omega$ and $y$ of $\KK^n$ such that $x + y \in \Omega$;

        \item\label{it:R-bounded} there exists a real number $M > 0$ such that
        \[ \abs{R(x, y)} \leqslant M \]
        for all elements $x$ of $\Omega$ and $y$ of $\KK^n$ such that $x + y \in \Omega$.
    \end{enumerate}
    Then there exist two integers $N$ and $\alpha$ with $\alpha \geqslant 1$, and some elements $a_1$, \dots, $a_n$ of~$\O_\KK^\times$ such that, for all elements $\lambda$ of $\KK$ with $\ord\lambda \leqslant N$, we have
    \[ I_{f, \varphi}(\lambda) = \psi(\lambda f(x_0))\varphi(x_0) \prod_{i = 1}^n \int_{\varpi_\KK^\alpha\O_\KK^{}} \psi(\lambda a_iu^2) \dd u. \]
    Moreover, the quantities $\alpha$ and $a_i$ are independent from the function $\varphi$.
\end{theo}

\begin{proof}
    We take a nonnegative integer $\alpha$ such that $B \coloneq x_0 + \varpi_\KK^\alpha\O_\KK^n \subset \Omega$. By Corollary~\ref{coro:Morse}, by assumptions~\ref{it:xO} and~\ref{it:xO-non-degen}, there exist a bianalytic function
    \[ T = (T_1, \dots, T_n) \colon (\varpi_\KK^\alpha\O_\KK^{})^n \to (\varpi_\KK^\alpha\O_\KK^{})^n \]
    with~$T(0) = 0$ and elements $a_1$, \dots, $a_n$ of $\O_\KK^\times$ such that, for all elements $x$ of $B$, we have
    \begin{equation}\label{eq:Morse-f}
        f(x) = f(x_0) + \sum_{i = 1}^n a_iT_i(x - x_0)^2.
    \end{equation}
    By additivity of integral, we cut the integral $I_{f, \varphi}(\lambda)$ in two parts according to the decomposition $\Omega = (\Omega \setminus B) \sqcup B$.

    On the open set $\Omega \setminus B$, the function $f$ has no critical point since it has only one critical point and the latter belongs to $B$. Thus, the function $x \longmapsto \abs{\grad f(x)}$ does not vanish on $\supp\varphi|_{\Omega \setminus B} = \supp\varphi \setminus B$. But this support is a closed subset of the compact set $\O_\KK^n$ and is so compact. Therefore, there exists a positive real number~$\delta$ such that
    \begin{equation}\label{eq:bounded-grad}
        \forall x \in \supp\varphi|_{\Omega \setminus B}, \qquad \abs{\grad f(x)} \geqslant \delta.
    \end{equation}
    By relation~\eqref{eq:bounded-grad} and by assumptions~\ref{it:def-R} and~\ref{it:R-bounded}, we use~\cite[Proposition 1.1]{Heifetz} (or equivalently~\cite[Proposition 2.5.3]{CHLR}): there exists an integer $N_1$ such that
    \[ \forall \lambda \in \KK^\times, \qquad \ord\lambda \leqslant N_1 \implies \int_{\Omega \setminus B} \varphi(x)\psi(\lambda f(x)) \dd x = 0. \]

    Now let us work on the open set $B$. Let $\lambda$ be an element of $\KK$. With decomposition~\eqref{eq:Morse-f}, the integral $B$ can be written as
    \[ \int_B \varphi(x)\psi(\lambda f(x)) \dd x = \psi(\lambda f(x_0)) \int_B \varphi(x) \prod_{i = 1}^n \psi(\lambda a_iT_i(x - x_0)^2) \dd x. \]
    The change of variables $y = T(x - x_0)$ gives
    \[ \int_B \varphi(x)\psi(\lambda f(x)) \dd x = \psi(\lambda f(x_0)) \int_{\KK^n} \theta(y)g_\lambda(y) \dd y \]
    where we consider the functions $\theta, g_\lambda \colon \KK^n \to \CC$ defined by the equality
    \[ \theta(y) \coloneq \begin{cases}
        \varphi(x_0 + T^{-1}(y))\abs{\det\Jac T^{-1}(y)} &\text{if } y \in (\varpi_\KK^\alpha\O_\KK^{})^n, \\
        0 &\text{else}
    \end{cases} \]
    and
    \[ g_\lambda(y) \coloneq \cf_{(\varpi_\KK^\alpha\O_\KK^{})^n}(y) \prod_{i = 1}^n \psi(\lambda a_iy_i^2) \]
    for all elements $y = (y_1, \dots, y_n)$ of $\KK^n$.

    The function $\theta$ is a Schwartz--Bruhat function because it has a compact support by definition and it is locally constant thanks to the following two points:
    \begin{itemize}
        \item since the function $\varphi$ is locally constant and the map $T$ is an open map, the function $y \longmapsto \varphi(x_0 + T^{-1}(y))$ is locally constant;

        \item the function $y \longmapsto \abs{\det\Jac T^{-1}(y)}$ is continuous and with discrete values (by formula~\eqref{eq:absolute-value}, its range is included in the discrete set $\{q_\KK^\ell\}_{\ell \in \ZZ}$ since the map~$T^{-1}$ is a diffeomorphism and so its jacobian is nonzero), so it is locally constant.
    \end{itemize}

    By Lemma~\ref{lemma:SB-Psi}, the function $g_\lambda$ is also a Schwartz--Bruhat function. Then using the inversion theorem (Theorem~\ref{thm:Fourier-inv}) and Corollary~\ref{coro:Plancherel}, we write
    \[ \int_{\KK^n} \theta(y)g_\lambda(y) \dd y = \int_{\KK^n} \hat{\hat\theta}(-y)g_\lambda(y) \dd y = \int_{\KK^n} \hat\theta(-\xi)\hat g_\lambda(\xi) \dd\xi. \]
    By Theorem~\ref{thm:Fourier-inv}, the function $\hat\theta$ is a Schwartz--Bruhat function and we consider an integer $\beta$ such that $\supp\hat\theta \subset (\varpi_\KK^\beta\O_\KK^{})^n$.

    Let $\xi = (\xi_1, \dots, \xi_n)$ be an element of $(\varpi_\KK^\beta\O_\KK^{})^n$. By Fubini's theorem, we write
    \begin{align*}
        \hat g_\lambda(\xi) &= \int_{\KK^n} \cf_{(\varpi_\KK^\alpha\O_\KK^{})^n}(y) \psi(\inner y\xi) \prod_{i = 1}^n \psi(\lambda a_iy_i^2) \dd y \\
        &= \prod_{i = 1}^n \int_{\varpi_\KK^\alpha\O_\KK^{}} \psi(\xi_iy_i + \lambda a_i^{}y_i^2) \dd y_i.
    \end{align*}
    Let $i$ be an index in $\{1, \dots, n\}$. We assume that $\lambda \ne 0$ and
    \begin{equation}\label{eq:bound-N2}
        \ord\lambda \leqslant N_{2, i} \coloneq \min\{\beta - \ord a_i - \alpha, 2\beta - \ord a_i - 1\}.
    \end{equation}
    Since $\ord\xi_i \geqslant \beta$, we have
    \[ \ord\xi_i - \ord\lambda a_i \geqslant \alpha \qquad \text{and} \qquad
       2\ord\xi_i - \ord\lambda a_i \geqslant 1. \]
    Thus, using Lemma~\ref{lemma:Gauss}, we get
    \[ \int_{\varpi_\KK^\alpha\O_\KK^{}} \psi(\xi_iu + \lambda a_iu^2) \dd u = \int_{\varpi_\KK^\alpha\O_\KK^{}} \psi(\lambda a_iu^2) \dd u. \]
    To conclude, using Fubini's theorem, we write
    \begin{align*}
        \int_B \varphi(x)\psi(\lambda f(x)) \dd x &= \psi(\lambda f(x_0)) \int_{\supp\hat\theta} \hat\theta(-\xi) \prod_{i = 1}^n \int_{\varpi_\KK^\alpha\O_\KK^{}} \psi(\lambda a_iu^2) \dd u \dd\xi \\
        &= \psi(\lambda f(x_0)) \int_{\KK^n} \hat\theta(-\xi) \dd\xi \times \prod_{i = 1}^n \int_{\varpi_\KK^\alpha\O_\KK^{}} \psi(\lambda a_iu^2) \dd u.
    \end{align*}
    Setting $\gamma \coloneq \abs{\det\Jac T^{-1}(0)}$, Theorem~\ref{thm:Fourier-inv} gives
    \[ \int_{K^n} \hat\theta(-\xi) \dd\xi = \hat{\hat\theta}(0) = \theta(0) = \gamma\varphi(x_0) \]
    and then
    \[ \int_B \varphi(x)\psi(\lambda f(x)) \dd x = \gamma\psi(\lambda f(x_0))\varphi(x_0) \prod_{i = 1}^n \int_{\varpi_\KK^\alpha\O_\KK^{}} \psi(\lambda a_iu^2) \dd u. \]
    To finish, we just have to take $N \coloneq \min\{N_1, N_{2, 1}, \dots, N_{2, n}\}$ and the proposition is proven. Finally, by the chain rule formula, since the function $T$ is bianalytic, we have~$\gamma = 0$.
\end{proof}

We deduce from Theorem~\ref{thm:FPS-p} the case of finitely many nondegenerate critical points.

\begin{coro}\label{coro:FPS-p}
    We use settings of Definition~\ref{def:osc-int-p} and make the following assumptions:
    \begin{enumerate}[label=(\roman*)]
        \item the function $f$ admits finitely many critical points $x_0$, \dots, $x_\ell$ in $\Omega$;

        \item the latter belong to the set $\supp\varphi$ and they are nondegenerate modulo~$\varpi_\KK$;

        \item there exists a function $R \colon \Omega \times \KK^n \to \KK$ such that
        \[ f(x + y) = f(x) + \inner{\grad f(x)}{y} + \inner{R(x, y)y}{y} \]
        for all elements $x$ of $\Omega$ and $y$ of $\KK^n$ such that $x + y \in \Omega$;

        \item there exists a real number $M > 0$ such that
        \[ \abs{R(x, y)} \leqslant M \]
        for all elements $x$ of $\Omega$ and $y$ of $\KK^n$ such that $x + y \in \Omega$.
    \end{enumerate}
    Then there exist two integers $N$ and $\alpha$ with $\alpha \geqslant 1$, and some elements $a_{j, i}$ of $\O_\KK^\times$ with~$0 \leqslant j \leqslant \ell$ and $1 \leqslant i \leqslant n$ such that, for all elements~$\lambda$ of $\KK$ with $\ord\lambda \leqslant N$, we have
    \[ I_{f, \varphi}(\lambda) = \sum_{j = 0}^\ell \psi(\lambda f(x_j))\varphi(x_j) \prod_{i = 1}^n \int_{\varpi_\KK^\alpha\O_\KK^{}} \psi(\lambda a_{j, i}u^2) \dd u. \]
    Moreover, the quantities $\alpha$ and $a_{j, i}$ are independent from the function $\varphi$.
\end{coro}

\begin{proof}
    We consider some disjoint balls around critical points and we use the additivity of integral: on balls, we use Theorem~\ref{thm:FPS-p} and, outside, we apply the stationary phase formula in the smooth case (see~\cite[Proposition~1.1]{Heifetz}).
\end{proof}

\begin{rema}
    The hypothesis $p_\KK \ne 2$ gets involved twice: one time in the proof of Theorem~\ref{thm:Morse} (it allows to symmetrize the problem and find some globally analytic square root function) and one another in the proof of Lemma~\ref{lemma:Gauss}. Without this hypothesis, these points seem to be compromised.
\end{rema}

\section{The motivic case}
\label{sec:motivic}

\subsection{A quick reminder about motivic integration}

In the two following paragraphs, we give a brief overview of the Cluckers--Loeser motivic's integration~\cite{CluckersLoeser2008}, its version with exponential~\cite{CluckersLoeser2010} and the Cluckers--Halupczok's notion of evaluation of constructible motivic functions~\cite{CluckersHalupczok}. We refer to these articles for more details.

\subsubsection{Cluckers--Loeser's motivic integration}

Let $k$ be a characteristic zero field. We consider a Denef-Pas language $\L_\text{DP}$, i.e. a language with three sorts: the valued field sort (equipped with the ring language), the residual field sort (equipped with an extension of the ring language by constant symbols) and the value group sort (equipped with an extension of the Presburger language by constant symbols) endowed with a symbol for the valuation $\ord$ and one for the angular map $\ac$. We consider the extension~$\L_\text{DP}(k)$ of this language by adding constant symbols for each element of the sets $k\ls t$ and $k$, respectively in the valued field sort and in the residue field sort.

\begin{exem}
    Let $K$ be a field. Endowed with the $t$-adic valuation
    \[ \ord \colon K\ls t \to \ZZ \cup \{+\infty\} \]
    and with the natural angular component
    \[ \ac \colon K\ls t \to K, \]
    the field $K\ls t$ is a valued field whose valuation ring is the ring $K\fs t$ of formal powers series and whose residual field is the field $K$. The triplet $(K\ls t, K, \ZZ)$ is a~$\L_\text{DP}$-structure.
\end{exem}

Let $\Field_k$ be the category of fields which contain the field $k$. Let $n$, $\ell$ and $r$ be three nonnegative integers. For a field $K$ of $\Field_k$, we set
\[ h[n, \ell, r](K) \coloneqq K\ls t^n \times K^\ell \times \ZZ^r. \]
A \emph{definable subassignment} $X$, also called \emph{definable set}, of $h[n, \ell, r]$ is a collection of subsets~$X(K)$ of $h[n, \ell, r](K)$ for each field $K$ of $\Field_k$ such that there exists an~$\L_\text{DP}(k)$-formula $\phi(x, \xi, \alpha)$ such that, for all fields $K$ of $\Field_k$, we have
\[ X(K) = \{(x, \xi, \alpha) \in h[n, \ell, r](K) \mid (K\ls t, K, \ZZ) \models \phi(x, \xi, a)\}. \]
A point of a definable subassignment $X$ is the data of a field $K$ of $\Field_k$ and an element of $X(K)$. We denote by $\abs X$ the set of points of $X$. For two definable subassignments~$X$ and $Y$, we write $X \subset Y$ if $X(K) \subset Y(K)$ for all fields $K$ of~$\Field_k$.

Let $\Def_k$ be the category of definable subassignments over the field $k$. Let us recall some notations of~\cite{CluckersLoeser2008, CluckersLoeser2010}. We denote by $\{*\}$ the final object of the category~$\Def_k$. For an object $X$ of $\Def_k$, we denote by
\begin{itemize}
    \item $X[n, \ell, r]$ the product $X \times h[n, \ell, r]$ for some nonnegative integers $n$, $\ell$ and~$r$;

    \item $\Def_X$ the category of objects of $\Def_k$ over $X$, i.e. of morphisms $Y \to X$ of~$\Def_k$ for some object $Y$ of $\Def_k$;

    \item $\RDef_X$ the category of objects $Y \to X$ of $\Def_X$ where $Y \subset X[0, \ell, 0]$ for some nonnegative integer $\ell$;

    \item $\RDefexp_X$ the category of objects $[Y \to X, \xi, g]$ for some object $Y \to X$ of~$\RDef_X$ and some morphisms $\xi \colon Y \to h[0, 1, 0]$ and $g \colon Y \to h[1, 0, 0]$ of~$\Def_k$. These objects will be denoted by $[Y \to X]\e(\xi)\E(g)$;

    \item $\K(\RDefexp_X)$ the Grothendieck ring of $\RDefexp_X$ (see~\cite[§3.1]{CluckersLoeser2010});

    \item $\Cons(X)\expe$ (respectively $\cons(X)\expe$) the set of exponential constructible motivic functions (respectively Functions) on $X$;

    \item $\LL$ the element $[X[0, 1, 0] \to X, 0, 0]$ of the ring $\K(\RDefexp_X)$;

    \item for an object $\Lambda$ of $\Def_k$ and an object $f \colon X \to \Lambda$ of $\Def_\Lambda$ whose all fibers have the same dimension, we denote by $\Inte_\Lambda(X)\expe$ or by $\Inte_f(X)\expe$ (respectively~$\inte_\Lambda\cons(X \to \Lambda)\expe$) the set of integrable exponential constructible motivic functions (respectively Functions) on~$X$ with parameter in $\Lambda$ as defined in~\cite[§7]{CluckersLoeser2010} and~\cite[Theorem 10.1.1]{CluckersLoeser2008}. If the definable subassignment $\Lambda$ is the final object $\{*\}$ and the definable morphism~$f$ is the final morphism, these sets will be simply denoted by $\Inte(X)\expe$ and $\inte\cons(X)\expe$.
\end{itemize}

For a morphism $f \colon X \to Y$ of $\Def_k$, the pullback and pushforward morphisms are respectively denoted by
\[ f^* \colon \Cons(Y)\expe \to \Cons(X)\expe \qquad \text{and} \qquad
   f_! \colon \inte\cons(X)\expe \to \inte\cons(Y)\expe \]
where the second morphism is given by the functor of motivic integration constructed in~\cite[Theorem 10.1.1]{CluckersLoeser2008} and~\cite[Theorem 4.1.1]{CluckersLoeser2010}. The latter morphism $f_!$ will be also denoted by $\mu_f$.

Moreover, let $\Lambda$ be a definable subassignment. For a morphism $f \colon X \to Y$ a morphism of $\Def_\Lambda$, the pushforward relatively to $\Lambda$ will be denoted by
\[ f_{!\Lambda} \colon \inte_\Lambda\cons(X \to \Lambda)\expe \to \inte\cons_\Lambda(Y \to \Lambda)\expe \]
which is given by the function of relative motivic integration (see~\cite[Theorem 14.1.1]{CluckersLoeser2008} and~\cite[Theorem 4.3.1]{CluckersLoeser2010}). Let $f \colon X \to \Lambda$ be an object of $\Def_\Lambda$. Its pushforward is the morphism of abelian groups
\[ f_{!\Lambda} \colon \inte_\Lambda\cons(X \to \Lambda)\expe \to \inte\cons_\Lambda(\Lambda \to \Lambda)\expe \cong \Cons(\Lambda)\expe \]
which will by denoted by $\mu_f$. Finally, for a function $\varphi$ of $\Inte_f(X)\expe$, we denote by~$\mu_f(\varphi)$ the function $\mu_f([\varphi])$ where the notation $[\varphi]$ stands for the class of the function $\varphi$ in~$\inte_\Lambda(X \to \Lambda)\expe$. In this way, we obtain a morphism of abelian groups
\[ \mu_f \colon \Inte_f(X)\expe \to \Cons(\Lambda)\expe \]
also denoted by $\mu_\Lambda$ or $f_!$.

\begin{rema}\label{rem:extension-zero}
    Let $X$ and $Y$ be two definable subassignments of $\Def_k$ such that~$X \subset Y$. We consider the canonical injection morphism $i \colon X \to Y$. Then the ring morphism~$i_!$ is the extension by zero on Presburger functions on $X$ to $Y$ and the base change from~$\K(\RDefexp_X)$ to $\K(\RDefexp_Y)$ induced by the composition.
\end{rema}

\subsubsection{Cluckers--Halupczok's evaluation}
\label{sssec:eval}

We recall the notion of evaluation of exponential constructible motivic functions as defined in~\cite{CluckersHalupczok}. We use the same notation, especially the class~$\S$ is here the collection of $\L_\text{DP}(k)$-structures of the form~$(K\ls t, K, \ZZ)$ for some field $K$ of $\Field_k$.

For a point $x = (x_0, K_0)$ of some definable subassignment $X$ ---~i.e. an $\S$-definable set, we define the extension $\L(x)$ of the language $\L \coloneq \L_\text{DP}(k)$ by adding constants symbols for the entries of the tuple $x_0$ and we define $\S(x)$ the class of~$\L(x)$-structures~$(K\ls t, K, \ZZ, x_{0K})$ for a field $K$ of $\Field_k$ which are elementary equivalent to the~$\L(x)$-structure $(K_0\ls t, K_0, \ZZ, x_0)$.

For a function~$f$ of $\Cons(X)\expe$ defined on a definable subassignment $X$ and a point~$x$ of~$X$, we consider the induced $\S(x)$-definable set $X_{\S(x)}$ and the induced~$\S(x)$-definable function~$f_{\S(x)}$ of $\Cons(X_{\S(x)})\expe$ and we set
\[ f(x) \coloneq \iota^*(f_{\S(x)}) \in \Cons(\{x\})\expe \]
where the $\S(x)$-function $\iota \colon \{x\} \to X_{\S(x)}$ is the canonical inclusion.

\begin{theo}[Cluckers--Halupczok {\cite[Theorem 1]{CluckersHalupczok}}]\label{thm:eval}
    Let $X$ be an $\S$-definable set. Let~$f$ and $g$ be two functions of $\Cons(X)\expe$. Then the following items are equivalent:
    \begin{enumerate}
        \item $f = g$;
        \item $f(x) = g(x)$ for all points $x$ of $X$.
    \end{enumerate}
\end{theo}

\begin{rema}
    This notion of evaluation can be slightly modified to suit different contexts (for example, see \cite[§2.6]{ForeyLoeserWyss}) and it leads to similar results such as Theorem~\ref{thm:eval}.
\end{rema}

The following result is a direct consequence of the construction of evaluation, but for the convenience of the reader, to show how things work, we give details.

\begin{lemm}\label{lemma:ext-zero}
    Let $X$ and $Y$ be two definable subassignments of $\Def_k$ such that~$X \subset Y$. We consider the canonical injection~$i \colon X \to Y$. Let~$\varphi$ be a function of $\Inte_i(X)\expe$ and~$i_!\varphi$ be the pushforward function which lies in $\Cons(Y)\expe$. Let $y$ be a point of~$Y$. Then
    \[ i_!\varphi(y) = \begin{cases}
        \varphi(y) &\text{if } y \in \abs X, \\
        0 &\text{else.}
    \end{cases} \]
\end{lemm}

\begin{proof}
    We can assume that the function $\varphi$ is of the form
    \[ \varphi = [p_Z \colon Z \to X, \xi, g] \otimes c\alpha\LL^\beta \in \Cons(X)\expe \]
    for some definable subassignment $Z$ of $X[0, n, 0]$ with an integer $n \geqslant 0$, some definable morphisms $\xi \colon Z \to h[0, 1, 0]$, $g \colon Z \to h[1, 0, 0]$ and $\alpha, \beta \colon X \to h[0, 0, 1]$ and some element $c$ of $\ZZ[\LL, \LL^{-1}, (1/(1 - \LL^{-i}))_{i \geqslant 1}]$. By Remark~\ref{rem:extension-zero}, we have
    \[ i_!\varphi = [i \circ p_Z \colon Z \to Y, \xi, g] \otimes ci_!(\alpha\LL^\beta) \]
    where the ring morphism $i_! \colon \Pres(X) \to \Pres(Y)$ is the extension by zero. By definition, we have
    \[ (i_!\varphi)_{\S(y)} = [Z_{\S(y)} \to Y_{\S(y)}, \xi_{\S(y)}, g_{\S(y)}] \otimes ci_{\S(y)!}(\alpha_{\S(y)}\LL^{\beta_{\S(y)}}) \]
    which gives
    \[ i_!\varphi(y) = [Z_{\S(y)} \otimes_{Y_{\S(y)}} \{y\} \to \{y\}, \xi_{\S(y)} \circ p, g_{\S(y)} \circ p] \otimes ci_{\S(y)!}(\alpha_{\S(y)}\LL^{\beta_{\S(y)}}) \circ i_y \]
    where the $\S(y)$-morphism $p \colon Z_{\S(y)} \otimes_{Y_{\S(y)}} \{y\} \to Z_{\S(y)}$ is the canonical projection and where the $\S(y)$-morphism $i_y \colon \{y\} \to Y_{\S(y)}$ is the canonical inclusion. Moreover, we have
    \[ \varphi(y) = [Z_{\S(y)} \otimes_{X_{\S(y)}} \{y\} \to \{y\}, \xi_{\S(y)} \circ p', g_{\S(y)} \circ p'] \otimes c(\alpha_{\S(y)}\LL^{\beta_{\S(y)}}) \circ i_y \]
    where the $\S(y)$-morphism $p' \colon Z_{\S(y)} \otimes_{X_{\S(y)}} \{y\} \to Z_{\S(y)}$ is the canonical projection. If $y \notin \abs X$, then $i_{\S(y)!}(\alpha_{\S(y)}\LL^{\beta_{\S(y)}}) \circ i_y = 0$ by Remark~\ref{rem:extension-zero} and thus $i_!\varphi(y) = 0$. We assume now that~$y \in \abs X$. On the one hand, we have
    \[ c \cdot i_{\S(y)!}(\alpha_{\S(y)}\LL^{\beta_{\S(y)}}) \circ i_y = c \cdot (\alpha_{\S(y)}\LL^{\beta_{\S(y)}}) \circ i_y. \]
    On the other hand, we have a canonical $\S(y)$-isomorphism
    \[ Z_{\S(y)} \otimes_{X_{\S(y)}} \{y\} \simeq Z_{\S(y)} \otimes_{Y_{\S(y)}} \{y\} \]
    which implies the equality
    \begin{multline*}
        [Z_{\S(y)} \otimes_{Y_{\S(y)}} \{y\} \to \{y\}, \xi_{\S(y)} \circ p, g_{\S(y)} \circ p] \\
        = [Z_{\S(y)} \otimes_{X_{\S(y)}} \{y\} \to \{y\}, \xi_{\S(y)} \circ p', g_{\S(y)} \circ p']
    \end{multline*}
    and conclude that $i_!\varphi(y) = \varphi(y)$.
\end{proof}

\subsection{Definable constants and open definable subassignments}

\subsubsection{Definable constants}
\label{sssec:constants}

In the following definition, we extend in a definable way and uniformly in the field, the choice of an element in a set.

\begin{defi}
    A \emph{definable constant} of a definable subassignment $X$ of~$\Def_k$ is a definable morphism $\{*\} \to X$ of $\Def_k$.
\end{defi}

Let $X$ be a definable subassignment of $\Def_k$ such that $X(K) \ne \emptyset$ for all fields~$K$ of~$\Field_k$. Let $c$ be a definable constant of $X$. For any field $K$ in $\Field_k$, we also denote by $c_K$ the unique image of the map $c_K \colon \{*\}(K) \to X(K)$. We denote by~$\{c\}$ the definable subassignment $(\{c_K\})_{K \in \Field_k}$.

\begin{rema}
    An element $c_0$ of the set $X(k)$ defines a definable constant $c$ of~$X$ by setting $c_K = c_0$ for any field $K$ of $\Field_k$. But all definable constants are not made in that way. Indeed, we give an example taking the field $k = \QQ$ and the definable constant $c$ of $\{0, 1\} \subset h[0, 0, 1]$ defined by the equivalence
    \[ c = 1 \iff \exists y \in h[1, 0, 0], \; y^2 = -1. \]
    Then $c_\QQ = 0 \ne c_\CC = 1$.
\end{rema}

\subsubsection{Open definable subassignments}

Let $\Lambda$ be a definable subassignment of $\Def_k$.

\begin{defi}
    Let $c \colon \Lambda \to h[1, 0, 0]$ and $\alpha \colon \Lambda \to h[0, 0, 1]$ be two definable morphisms of $\Def_k$. The \emph{ball} of valuative radius $\alpha$ and centre $c$ is the definable subassignment $\B_\Lambda(c, \alpha)$ of $\Lambda[1, 0, 0]$ defined, for any field $K$ of $\Field_k$, by
    \[ \B_\Lambda(c, \alpha)_K \coloneq \{(\lambda, x) \in \Lambda(K) \times K\ls t \mid \ord(x - c_K(\lambda)) \geqslant \alpha_K(\lambda)\}. \]
    Endowed with the canonical projection $\B_\Lambda(c, \alpha) \to \Lambda$, it is an object of $\Def_\Lambda$.

    When the definable subassignment $\Lambda$ is the final object $\{*\}$, this ball will be simply denoted by $\B(c, \alpha)$. We denote by $\O$ the ball $\B(0, 0) = (K\fs t)_{K \in \Field_k}$.
\end{defi}

\begin{defi}\label{def:def-open}
    A definable subassignment $\Omega$ of $\Lambda[n, 0, 0]$ is \emph{$\Lambda$-open} if, for all points~$(\lambda, K)$ of $\Lambda$, the fiber
    \[ \Omega(K)_\lambda \coloneq \{x \in K\ls t^n \mid (\lambda, x) \in \Omega(K)\} \subset K\ls t^n \]
    is open in $K\ls t^n$.

    When the definable subassignment $\Lambda$ is the final object $\{*\}$, we say that the definable subassignment $\Omega$ is \emph{open}.
\end{defi}

\begin{rema}
    Let $K$ a field of $\Field_k$. In Definition~\ref{def:def-open}, we do not specify the topology on the space $K\ls t^n$ because, since the field $K\ls t^n$ is complete for the $t$-adic topology, all norms on the finite-dimensional vector space $K\ls t^n$ over the field~$K\ls t$ are equivalent. We will use the norm
    \[ (x_1, \dots, x_n) \in K\ls t^n \longmapsto \max(\operatorname{exp}(-\ord x_1), \dots, \operatorname{exp}(-\ord x_n)) \]
    on the vector space $K\ls t^n$.
\end{rema}

\begin{prop}\label{prop:def-open}
    Let $\Omega$ be a $\Lambda$-open definable subassignment of $\Lambda[n, 0, 0]$ with
    \[ \text{for all fields $K$ in $\Field_k$}, \qquad \Omega(K) \ne \emptyset \]
    and $x \colon \Lambda \to h[n, 0, 0]$ be a morphism of $\Def_k$ whose components are $x_1$, \dots, $x_n$. We assume that
    \[ \text{for all points $(\lambda, K)$ of $\Lambda$}, \qquad (\lambda, x_K(\lambda)) \in \Omega(K). \]
    Then there exists a definable morphism $\alpha \colon \Lambda \to h[0, 0, 1]$ of $\Def_k$ such that
    \[ \B_\Lambda(x_1, \alpha) \times_\Lambda \dots \times_\Lambda \B_\Lambda(x_n, \alpha) \subset \Omega. \]
\end{prop}

\begin{proof}
    For each point $(\lambda, K)$ of $\Lambda$, we set
    \[ \alpha_K(\lambda) \coloneq \max\{0, \inf\{\beta \in \ZZ \mid x_K(\lambda) + t^\beta K\ls t^n \subset \Omega(K)_\lambda\}\} \]
    which is finite since the set $\Omega(K)_\lambda$ is open and we have $x_K(\lambda) \in \Omega(K)_\lambda$. Thus, we get a convenient definable morphism $\alpha \colon \Lambda \to h[0, 0, 1]$.
\end{proof}

\subsection{Globally analytic functions and critical points}

We define the definable subassignment $\M_n(\O)$ of $\Def_k$ as the object~$\O^{n^2}$ of $\Def_k$ where, for any field~$K$ of $\Field_k$, we consider the set $\M_n(\O)(K)$ as the set $\M_n(K\fs t)$ of square matrices of dimension~$n$ with coefficients in $K\fs t$. Using Remark~\ref{rem:gen-Kt}, we give the following definition.

\begin{defi}\label{def:analytic-mot1}
    Let $A$ be an open definable subassignment of $\O^n$ and $B$ be an open definable subassignment of~$\O^m$.
    \begin{itemize}
        \item A definable morphism $A \to \O$ is \emph{globally analytic} if, for any field $K$ of~$\Field_k$, the function $f_K \colon A(K) \to K\fs t$ is globally analytic.

        \item A definable morphism $A \to \O^m$ is globally analytic if its components are globally analytic.

        \item A definable morphism $A \to B$ is a \emph{bianalytic definable morphism} if it is a globally analytic definable isomorphism whose inverse is globally analytic.
    \end{itemize}
\end{defi}

As we can define partial derivatives of a function with a first order logical formula, we can also give the following definition.

\begin{defi}\label{def:analytic-mot2}
    Let $f \colon A \to \O$ be a globally analytic definable morphism and $\alpha$ be a definable constant of $h[0, 0, 1]$ such that $\alpha_K \geqslant 1$ for all fields $K$ of~$\Field_k$.
    \begin{itemize}
        \item The \emph{gradient} at some definable constant $x$ of $A$ of the morphism $f$ is the definable constant $\grad f(x)$ of $\O^n$ given by the equality
        \[ \forall K \in \Field_k, \qquad [\grad f(x)]_K \coloneq \grad f_K(x_K) \]
        and its \emph{Hessian} at some definable constant $x$ of $A$ is the definable constant~$\Hess f(x)$ of $\M_n(\O)$ given by the equality
        \[ \forall K \in \Field_k, \qquad [\Hess f(x)]_K \coloneq \Hess f_K(x_K). \]

        \item A \emph{critical point} of the morphism $f$ is a definable constant $x$ of $A$ such that
        \[ \forall K \in \Field_k, \qquad [\grad f(x)]_K = 0. \]

        \item A critical point $x$ of the morphism $f$ is \emph{nondegenerate} if
        \[ \forall K \in \Field_k, \qquad \det([\Hess f(x)]_K) \ne 0. \]

        \item A \emph{critical point modulo $t^\alpha$} of the morphism $f$ is a definable constant $x$ of~$A$ such that
        \[ \forall K \in \Field_k, \qquad [\grad f(x)]_K \equiv 0 \mod{t^{\alpha_K}}. \]

        \item A critical point $x$ of the morphism $f$ is \emph{nondegenerate modulo $t$} if
        \[ \forall K \in \Field_k, \qquad \det([\Hess f(x)]_K) \not\equiv 0 \mod t. \]
    \end{itemize}
\end{defi}

\subsection{A motivic Morse lemma}

Let $K$ be a field of $\Field_k$.

\begin{prop}\label{prop:Morse-Laurent}
    Let $\Omega$ be an open set of $K\fs t^n$. Let $f \colon \Omega \to K\fs t$ be a globally analytic function with a critical point $x_0$ in $\Omega$ which is nondegenerate modulo $t$. Let~$\alpha$ be a positive integer such that $B \coloneq x_0 + t^\alpha K\ls t^n \subset \Omega$. Then the point $x_0$ is the only critical point of the function $f$ on the polydisc $B$. Moreover, there exist a bianalytic function
    \[ T = (T_1, \dots, T_n) \colon (t^\alpha K\fs t)^n \to (t^\alpha K\fs t)^n \]
    with $T(0) = 0$ and elements $a_1$, \dots, $a_n$ of $K\fs t^\times$ such that, for all elements $x$ of $B$, we have
    \[ f(x) = f(x_0) + \sum_{i = 1}^n a_iT_i(x - x_0)^2. \]
\end{prop}

\begin{proof}
    The result follows from an analogue of Theorem~\ref{thm:Morse} (and Corollary~\ref{coro:Morse}) with the field $\KK = K\ls t$. The proof is the same since the characteristic of the field~$K$ is not equal to $2$.
\end{proof}

We can now state a motivic Morse lemma which is the definable version of Proposition~\ref{prop:Morse-Laurent}.

\begin{prop}\label{prop:Morse-mot}
    Let $\alpha$ be a definable constant of $h[0, 0, 1]$ such that
    \[ \text{for all fields $K$ of $\Field_k$,} \qquad \alpha_K \geqslant 1. \]
    Let $f \colon \B(0, \alpha)^n \to \O$ be a globally analytic definable morphism such that the origin is a critical point modulo~$t^\alpha$ which is nondegenerate modulo $t$. Then the morphism~$f$ admits a unique critical point $c$ in $\B(0, \alpha)^n$ and there exist a bianalytic definable morphism
    \[ T = (T_1, \dots, T_n) \colon \B(0, \alpha)^n \to \B(0, \alpha)^n \]
    with $T(0) = 0$ and definable constants $a_1$, \dots, $a_n$ of $\O^\times$ such that, for all points $(x, K)$ of $\B(0, \alpha)^n$, we have
    \[ f_K(x) = f_K(c_K) + \sum_{i = 1}^n a_{iK}T_{iK}(x)^2. \]
\end{prop}

\begin{proof}
    Applying Proposition~\ref{prop:Morse-Laurent}, for each $K$ in $\Field_k$, we get a critical point $c_K$, a bianalytic map $T_K$ and constants $a_{iK}$. By definability of the gradient and the constructions in the proof of Proposition~\ref{prop:Morse-Laurent}, we get a collection $c$ which is a definable constant, a collection $T$ which is a bianalytic definable morphism and collections $a_i$ which are definable constants.
\end{proof}

\begin{coro}\label{coro:Morse-mot}
    Let $\Omega$ be an open definable subassignment of $\O^n$ and let~$f \colon \Omega \to \O$ be a globally analytic definable morphism with some critical point~$x_0 = (x_1, \dots, x_n)$ in $\Omega$ which is nondegenerate modulo $t$. Let $\alpha \geqslant 1$ be a definable constant of~$h[0, 0, 1]$ such that $B \coloneq \B(x_1, \alpha) \times \dots \times \B(x_n, \alpha) \subset \Omega$. Then the definable constant $x_0$ is the only critical point of the morphism $f$ on the polydisc~$B$. Moreover, there exist a bianalytic definable morphism
    \[ T = (T_1, \dots, T_n) \colon \B(0, \alpha)^n \to \B(0, \alpha)^n \]
    with $T(0) = 0$ and definable constants $a_1$, \dots, $a_n$ of $\O^\times$ such that, for all points $(x, K)$ of~$\B(0, \alpha)^n$, we have
    \[ f_K(x) = f_K(x_{0K}) + \sum_{i = 1}^n a_{iK}T_{iK}(x - x_{0K})^2. \]
\end{coro}

\begin{proof}
    The existence of the definable constant $\alpha$ follows from Proposition~\ref{prop:def-open}. Then translating, the result follows from Proposition~\ref{prop:Morse-mot}.
\end{proof}

\subsection{Schwartz--Bruhat functions on an open definable subassignment}

We first recall the notion of Schwartz--Bruhat functions from~\cite[§7.5]{CluckersLoeser2010}. Let $\Lambda$ be a definable subassignment.

\begin{defi}\label{def:SB}
    A function $\varphi$ of $\Inte(h[n, 0, 0])\expe$ is a \emph{Schwartz--Bruhat function} on~$h[n, 0, 0]$ if the two following points are satisfied:
    \begin{enumerate}[label=(\roman*)]
        \item\label{it:bounded-support} there exists a definable morphism $\alpha^-(\varphi) \colon \{*\} \to h[0, 0, 1]$ such that
        \[ \varphi \cdot \cf_{\B_\Lambda(0, \alpha)^n} = \varphi \]
        for all definable morphisms $\alpha \colon \{*\} \to h[0, 0, 1]$ such that $\alpha \leqslant \alpha^-(\varphi)$;

        \item\label{it:local-constancy} there exists a definable morphism $\alpha^+(\varphi) \colon \{*\} \to h[0, 0, 1]$ such that
        \[ \varphi * \cf_{\B_\Lambda(0, \alpha)^n} = \LL^{-n\alpha}\varphi \]
        for all definable morphisms $\alpha \colon \{*\} \to h[0, 0, 1]$ such that $\alpha \geqslant \alpha^+(\varphi)$ where the symbol $*$ denotes the convolution product (see~\cite[§7.4]{CluckersLoeser2010}).
    \end{enumerate}
\end{defi}

\begin{rema}
    In Definition~\ref{def:SB}, informally, point~\ref{it:bounded-support} ensures that the support of the function $\varphi$ is bounded and point~\ref{it:local-constancy} ensures that the function $\varphi$ is locally constant. This corresponds closely to the $p$-adic case described above.
\end{rema}

Let $\Omega$ be an open definable subassignment of $h[n, 0, 0]$. We consider the canonical injection $i \colon \Omega \to h[n, 0, 0]$.

\begin{defi}
    A function $\varphi$ of $\Inte(\Omega)\expe$ is a \emph{Schwartz--Bruhat function} on~$\Omega$ if the function $i_!\varphi$\footnote{A function of $\Inte(\Omega)\expe$ is also a function of $\Inte_i(\Omega)\expe$ by~\cite[Theorem 10.1.1, Axiom A0.(b)]{CluckersLoeser2008}, so we can take the pushforward function.} is a Schwartz--Bruhat function on $h[n, 0, 0]$.
\end{defi}

In this case, we set $\alpha^-(\varphi) \coloneqq \alpha^-(i_!\varphi)$ and $\alpha^+(\varphi) \coloneqq \alpha^+(i_!\varphi)$. We denote by~$\Sch(\Omega)\expe$ the set of Schwartz--Bruhat functions on $\Omega$.

In the following proposition (see also~\cite{RaibautForey}), we consider the evaluation point of view of point~\ref{it:local-constancy} of Definition~\ref{def:SB} on Schwartz--Bruhat functions.

\begin{prop}\label{prop:SB-eval}
    Let $\varphi $ be a function of $\Inte(h[n, 0, 0])\expe$ and $\alpha$ be a definable constant of $h[0, 0, 1]$. We set
    \[ X \coloneq \{(x, z) \in h[2n, 0, 0] \mid \ord(x - z) \geqslant \alpha\} \]
    and we consider the two projections $p_x, p_z \colon X \to h[n, 0, 0]$. Then the following items are equivalent:
    \begin{enumerate}[label=(\roman*)]
        \item\label{it:convol} $\varphi * \cf_{\B(0, \alpha)^n} = \LL^{-\alpha n}\varphi$;

        \item\label{it:proj} $p_x^*\varphi = p_z^*\varphi$;

        \item\label{it:eval} for all points $(x, z, K)$ of $X$, we have $\varphi(x, K)_{\S(x, z, K)} = \varphi(z, K)_{\S(x, z, K)}$.
    \end{enumerate}
\end{prop}

In Proposition~\ref{prop:SB-eval}, we have identified the three classes $\S(x, z, K)$, $\S(z, K)(x, K)$ and $\S(x, K)(z, K)$.

\begin{proof}
    Let us prove the equivalence $\text{\ref{it:proj}} \iff \text{\ref{it:eval}}$. By~\cite[Corollary 3.6.1]{CluckersHalupczok}, for all points $(x, z, K)$ of $X$, we can write
    \[ p_x^*\varphi(x, z, K) = \varphi(p_x(x, z, K))_{\S(x, z, K)} = \varphi(x, K)_{\S(x, z, K)} \]
    and
    \[ p_z^*\varphi(x, z, K) = \varphi(p_z(x, z, K))_{\S(x, z, K)} = \varphi(z, K)_{\S(x, z, K)}. \]
    Then the equivalence follows from Theorem~\ref{thm:eval}

    Let us prove the implication $\text{\ref{it:proj}} \implies \text{\ref{it:convol}}$. We assume that $p_x^*\varphi = p_z^*\varphi$. For each point $(x, z, K)$ of $X$, we can write
    \[ \cf_{\B(0, \alpha)^n}(x - z, K) = \cf_X(x, z, K). \]
    By definition of the convolution product, we thus get
    \[ \varphi * \cf_{\B(0, \alpha)^n} = \mu_z(p_x^*\varphi \cdot \cf_X). \]
    Then the hypothesis $p_x^*\varphi = p_z^*\varphi$ implies
    \[ \varphi * \cf_{\B(0, \alpha)^n} = \mu_z(p_z^*\varphi \cdot \cf_X). \]
    Finally, the projection axiom (A3) of~\cite[Theorem 4.1.1]{CluckersLoeser2010} implies
    \[ \varphi * \cf_{\B(0, \alpha)^n} = \varphi \cdot \mu_z(\cf_X) = \LL^{-\alpha n}\varphi. \]

    Let us prove the implication $\text{\ref{it:convol}} \implies \text{\ref{it:eval}}$. Let $(x, z, K)$ be a point of $X$. As we have~$\ord(x - z) \geqslant \alpha$, the $\S(x, z, K)$-definable sets $\B(x, \alpha_{\S(x, z, K)})$ and~$\B(z, \alpha_{\S(x, z, K)})$ are equal and, applying~\cite[Corollary 3.6.1]{CluckersHalupczok}, we can write
    \begin{align*}
        \varphi * \cf_{\B(0, \alpha)^n}(z, K)_{\S(x, z, K)} &= \left(\int_{p_z^{-1}(z)} (p_x^*\varphi \cdot \cf_X)|_{p_z^{-1}(z)}\right)_{\S(x, z, K)} \\
        &= \left(\int_{h[n, 0, 0]} \varphi(y)\cf_X(y, z) \dd y\right)_{\S(x, z, K)} \\
        &= \left(\int_{h[n, 0, 0]} \varphi(x)\cf_{\B(z, \alpha)}(y) \dd y\right)_{\S(x, z, K)} \\
        &= \left(\int_{h[n, 0, 0]} \varphi(x)\cf_{\B(x, \alpha)}(y) \dd y\right)_{\S(x, z, K)} \\
        &= \varphi * \cf_{\B(0, \alpha)^n}(x, K)_{\S(x, z, K)}.
    \end{align*}
    With assumption~\ref{it:convol}, we deduce that $\LL^{-\alpha n}\varphi(z, K)_{\S(x, z, K)} = \LL^{-\alpha n}\varphi(x, K)_{\S(x, z, K)}$ and so $\varphi(z, K)_{\S(x, z, K)} = \varphi(x, K)_{\S(x, z, K)}$ by multiplying by $\LL^{\alpha n}$.
\end{proof}

\begin{defi}
    Let $\varphi$ be a function of $\Inte(h[n, 0, 0])\expe$. The function $\varphi$ is said to be \emph{locally constant} if there exists a definable constant $\alpha$ of~$h[0, 0, 1]$ such that the equivalent conditions of Proposition~\ref{prop:SB-eval} are satisfied.
\end{defi}

\begin{prop}\label{prop:prod-SB}
    The product of two Schwartz--Bruhat functions on an open definable subassignment of $h[n, 0, 0]$ is a Schwartz--Bruhat function.
\end{prop}

\begin{proof}
    We denote $i \colon \Omega \to h[n, 0, 0]$ the inclusion morphism. We use the fact that the map $i_!$ is a ring morphism and we apply~\cite[Proposition 2.27]{Raibaut}.
\end{proof}

\subsection{The formula}
\label{sec:formula-mot}

Following~\cite[§4]{Raibaut}, we set up notation for the oscillatory integral we will work on.

\begin{defi}\label{def:osc-int-mot}
    Let $\Omega$ be an open definable subassignment of~$\O^n$. Let~$\varphi$ be a Schwartz--Bruhat function on $\Omega$ and $f \colon \Omega \to \O$ be a globally analytic definable morphism. Let~$\Lambda$ be an unbounded definable subassignment of~$h[1, 0, 0]$. We consider the definable morphism
    \[ \fonc{\lambda f}{\Lambda \times \Omega}{h[1, 0, 0],}{(\lambda, x)}{\lambda f(x).} \]
    Let $\pi \colon \Lambda \times \Omega \to \Lambda$ and $\pi' \colon \Lambda \times \Omega \to \Omega$ be the canonical projections. Finally, we consider the function
    \[ I_{f, \varphi} \coloneq \pi_!(\pi'^*(\varphi)\E(\lambda f)) \in \Cons(\Lambda)\expe. \]
\end{defi}

\begin{rema}
    Until the end, we will use the integral notation. In particular, thanks to evaluation in §\ref{sssec:eval}, for all points $\lambda$ of $\Lambda$, we can write
    \[ I_{f, \varphi}(\lambda) = \int_\Omega \varphi(x)\E(\lambda f(x)) \dd x \]
    in the ring $\Cons(\{\lambda\})\expe$.
\end{rema}

\begin{theo}\label{thm:FPS-mot}
    We use the setting and notation of Definition~\ref{def:osc-int-mot} and make the following assumptions:
    \begin{enumerate}[label=(\roman*)]
        \item\label{it:mot-x0K} for any field $K$ of $\Field_k$, the function $f_K$ has a unique critical point $x_{0K}$ in~$\Omega(K)$, this critical point is nondegenerate modulo $t$ and belongs to the polydisc~$\B(0, \alpha^-(\varphi))^n(K)$;

        \item\label{it:mot-x0} the family $(x_{0K})_{K \in \Field_k}$ defines a definable constant of $\Omega$, denoted by $x_0$;

        \item\label{it:mot-def-R} there exists a definable morphism $R \colon \Omega[n, 0, 0] \to h[1, 0, 0]$ such that, for all fields $K$ of $\Field_k$, we have
        \[ f_K(x + y) = f_K(x) + \inner{\grad f_K(x)}{y} + \inner{R_K(x, y)y}{y} \]
        for all elements $x$ of $\Omega(K)$ and $y$ of $K\ls t^n$ with $x + y \in \Omega(K)$;

        \item\label{it:mot-R-bounded} there exists a definable constant $M$ of $h[0, 0, 1]$ such that, for all fields $K$ of~$\Field_k$, we have
        \[ \ord R_K(x, y) \geqslant M_K \]
        for all elements $x = (x_1, \dots, x_n)$ of $\Omega(K)$ and $y = (y_1, \dots, y_n)$ of $K\ls t^n$ with~$\ord x_i \geqslant \alpha^-(\varphi)_K$ and $\ord y_i \geqslant \alpha^+(\varphi)_K$ for all indexes $i$ in $\{1, \dots, n\}$.
    \end{enumerate}
    Then there exist definable constants $A$ and $\alpha$ of $h[0, 0, 1]$ with $\alpha \geqslant 1$ and definable constants $a_1$, \dots, $a_n$ of $\O^\times$ such that, in the ring $\Cons(\Lambda)\expe$, we have
    \[ \cf_\Gamma \int_\Omega \varphi(x)\E(\lambda f(x)) \dd x = \cf_\Gamma \E(\lambda f \circ x_0)i_{x_0}^*(\varphi) \prod_{i = 1}^n \left(\int_{\B(0, \alpha)} \E((\lambda, y_i) \longmapsto \lambda a_iy_i^2) \dd y_i\right) \]
    with $\Gamma \coloneq \{\lambda \in \Lambda \mid \ord\lambda \leqslant A\}$ and where the morphism $i_{x_0} \colon \{x_0\} \to \Omega$ is the canonical inclusion morphism of the definable subassignment $\{x_0\}$ (see §\ref{sssec:constants}).
\end{theo}

\begin{rema}
    As in the non-Archimedean local field case (see Corollary~\ref{coro:FPS-p}), we get immediately a corollary when the morphism $f$ has finitely many nondegenerate modulo~$t$ critical points. Furthermore, Theorem~\ref{thm:FPS-mot} has also a version with parameters, i.e. for Schwartz--Bruhat functions $\varphi$ of $\S_P(\Omega)$ with definable subassignments~$P$ and~$\Omega \subset P[n, 0, 0]$ (see~\cite[§7.5]{CluckersLoeser2010}).
\end{rema}

\begin{proof}[Proof of Theorem~\ref{thm:FPS-mot}]
    In the proof, we will use evaluation as §\ref{sssec:eval} and integral notation. We denote~$x_0 = (x_1, \dots, x_n)$. Let~$\alpha$ be a definable constant of~$h[0, 0, 1]$ which is positive such that
    \[ B \coloneq \B(x_1, \alpha) \times \dots \times \B(x_n, \alpha) \subset \Omega. \]
    By Corollary~\ref{coro:Morse-mot}, thanks to assumptions~\ref{it:mot-x0K} and~\ref{it:mot-x0}, there exist a bianalytic definable morphism
    \[ T = (T_1, \dots, T_n) \colon \B(0, \alpha)^n \to \B(0, \alpha)^n \]
    with $T(0) = 0$ and definable constants $a_1$, \dots, $a_n$ of $\O^\times$ such that
    \begin{equation}\label{eq:dec}
        f_K(x) = f_K(x_{0K}) + \sum_{i = 1}^n a_{iK}T_{iK}(x - x_{0K})^2
    \end{equation}
    for all points $(x, K)$ of $\B(0, \alpha)^n$.

    Let $\lambda$ a point of $\Lambda$. We work now over the residual field of the point $\lambda$. We split the integral $I_{f, \varphi}(\lambda)$ in two pieces: we write
    \[ I_{f, \varphi}(\lambda) = \int_\Omega \cf_{\Omega \setminus B}(x)\varphi(x)\E(\lambda f(x)) \dd x + \int_\Omega \cf_B(x)\varphi(x)\E(\lambda f(x)) \dd x. \]

    \smallbreak

    We work on the first integral. First, the function $\cf_{\Omega \setminus B}\varphi$ is a Schwartz--Bruhat function since it is a product of two Schwartz--Bruhat functions by Proposition~\ref{prop:SB-eval}.

    Let us prove that there exists a definable constant $N$ of $h[0, 0, 1]$ such that, for all fields $K$ of~$\Field_k$, we have
    \[ \ord\grad f_K(x) \leqslant N \]
    for all elements $x = (x_1, \dots, x_n)$ of $\Omega(K)$ with $\ord x_i \geqslant \alpha^-(\varphi)_K$ for all indexes $i$ in~$\{1, \dots, n\}$. By assumption~\ref{it:mot-x0K}, the definable morphism $\grad f$ does not vanish on~$\Omega \setminus \B(0, \alpha^-(\varphi))^n$. Then we get a definable continuous morphism
    \[ \ord\grad f \colon \Omega \setminus \B(0, \alpha^-(\varphi))^n \to h[0, 0, 1] \]
    on the closed and bounded definable subassignment $\Omega \setminus \B(0, \alpha^-(\varphi))^n$. Then applying~\cite[Proposition 3.3]{Raibaut}, we find a such definable constant $N$. Thus, by assumptions~\ref{it:mot-def-R} and~\ref{it:mot-R-bounded}, by~\cite[Proposition 4.2]{Raibaut}\footnote{Proposition 4.2 of~\cite{Raibaut} is stated when the definable subassignment $\Lambda$ is a cone
    \[ \Lambda_{\mathfrak n} \coloneqq \{x \in h[1, 0, 0] \mid \ord x \equiv 0 \mod{\mathfrak n}, \; \ac x = 1\} \]
    for some nonnegative integer $\mathfrak n \geqslant 1$, but it remains valid for general unbounded definable subassignments $\Lambda$ of~$h[1, 0, 0]$.}, there exists a definable constant $A_1$ of~$h[0, 0, 1]$ such that
    \[ \cf_{\Gamma_1}(\lambda) \int_\Omega \cf_{\Omega \setminus B}(x)\varphi(x)\E(\lambda f(x)) \dd x = 0 \]
    with $\Gamma_1 \coloneq \{\lambda \in \Lambda \mid \ord\lambda \leqslant -A_1 - N\}$.

    \smallbreak

    We work now on the second integral. By decomposition~\eqref{eq:dec}, we have
    \[ \int_\Omega \cf_B(x)\varphi(x)\E(\lambda f(x)) \dd x = \E(\lambda f \circ x_0) \int_B \varphi(x) \prod_{i = 1}^n \E(\lambda a_iT_i(x - x_0)^2) \dd x. \]
    We consider the definable morphism
    \[ \fonc{u}{\B(0, \alpha)^n}{B,}{y}{x_0 + T^{-1}(y).} \]
    By~\cite[Theorem 4.2.1]{CluckersLoeser2010}, with the change of variables $u$, we have
    \[ \int_\Omega \cf_B(x)\varphi(x)\E(\lambda f(x)) \dd x = \E(\lambda f \circ x_0) \int_{h[n, 0, 0]} \theta(y)g_\lambda(y) \dd y \]
    where we consider the functions
    \begin{equation}\label{eq:theta}
        \theta \coloneq \iota_!((j \circ u)^*(\varphi)\LL^{\ord\Jac u^{-1} \circ u}) \in \Cons(h[n, 0, 0])\expe
    \end{equation}
    and
    \[ \Inte(h[n, 0, 0])\expe \ni g_\lambda \colon y = (y_1, \dots, y_n) \longmapsto \cf_{\B(0, \alpha)^n}(y) \prod_{i = 1}^n \E(\lambda a_iy_i^2) \]
    where the morphisms $\iota \colon \B(0, \alpha)^n \to h[n, 0, 0]$ and $j \colon B \to \Omega$ are the inclusion morphisms. By Proposition~\ref{prop:theta-SB}, the function $\theta$ is a Schwartz--Bruhat function and, by Theorem~\ref{prop:SB-eval} and Axiom (R4) for~\cite{CluckersLoeser2010}, the function $g_\lambda$ is a Schwartz--Bruhat function. By an analogue result of Corollary~\ref{coro:Plancherel} and by~\cite[Theorem 7.5.1]{CluckersLoeser2010}, we have
    \[ \int_{h[n, 0, 0]} \theta(y)g_\lambda(y) \dd u = \int_{h[n, 0, 0]} \LL^n\hat{\hat\theta}(-y)g_\lambda(y) \dd y = \LL^n \int_{h[n, 0, 0]} \hat\theta(-\xi)\hat g_\lambda(\xi) \dd\xi \]
    where the hat stands for the Fourier transform (see~\cite[§7.2]{CluckersLoeser2010}).

    We study now the function $\hat g_\lambda$. By Fubini's theorem, i.e. the functoriality of the pushforward operation (see~\cite[Theorem 4.1.1]{CluckersLoeser2010}), we have
    \begin{align*}
        \cf_{\B(0, \alpha^-(\hat\theta))^n}(\xi)\hat g_\lambda(\xi) &= \cf_{\B(0, \alpha^-(\hat\theta))^n}(\xi) \int_{h[n, 0, 0]} \cf_{\B(0, \alpha)^n}(y)\E(\inner y\xi) \prod_{i = 1}^n \E(\lambda a_iy_i^2) \dd y \\
        &= \cf_{\B(0, \alpha^-(\hat\theta))^n}(\xi) \prod_{i = 1}^n \int_{\B(0, \alpha)} \E(\xi_iy_i + \lambda a_iy_i^2) \dd y_i.
    \end{align*}
    As in inequalities~\eqref{eq:bound-N2} in the proof of Theorem~\ref{thm:FPS-p}, using Lemma~\ref{lemma:Gauss-mot}, we can find a definable constant~$A_2$ of $h[0, 0, 1]$ such that
    \begin{multline*}
        \cf_{\B(0, \alpha^-(\hat\theta))}(\xi_i)\cf_{\Gamma_2}(\lambda) \int_{\B(0, \alpha)} \E(\xi_iy_i + \lambda a_iy_i^2) \dd y_i \\
        = \cf_{\B(0, \alpha^-(\hat\theta))}(\xi_i)\cf_{\Gamma_2}(\lambda) \int_{\B(0, \alpha)} \E(\lambda a_iy_i^2) \dd y_i
    \end{multline*}
    with $\Gamma_2 \coloneq \{\lambda \in \Lambda \mid \ord\lambda \leqslant A_2\}$ for all indexes $i$ in $\{1, \dots, n\}$. Finally, by Axiom~(A3) of~\cite[Theorem 4.1.1]{CluckersLoeser2010}, we find
    \[ \cf_{\Gamma_2}(\lambda) \int_B \varphi(x)\E(\lambda f(x)) \dd x = \cf_{\Gamma_2}(\lambda) \left(\int_{h[n, 0, 0]} \hat\theta(-\xi) \dd\xi\right) \prod_{i = 1}^n \left(\int_{\B(0, \alpha)} \E(\lambda a_iy_i^2) \dd y_i\right)\kern-\nulldelimiterspace. \]
    By analog arguments as in the proof of Proposition~\ref{prop:theta-SB}, we have
    \[ \ord\Jac T^{-1}(0) = 0. \]
    Then by~\cite[Theorem 7.5.1]{CluckersLoeser2008} and~\cite[§4.3]{CluckersLoeser2010} with $\Lambda = \{*\}$, we can write
    \[ \int_{h[n, 0, 0]} \hat\theta(-\xi) \dd\xi = i_0^*(\hat{\hat\theta}) = i_0^*(\LL^{-n}\theta) = \LL^{-n}i_{x_0}^*(\varphi) \]
    where the morphism $i_0 \colon \{0\} \to h[n, 0, 0]$ stands for the canonical injection. Thus, we obtain
    \[ \cf_{\Gamma_2}(\lambda) \int_B \varphi(x)\E(\lambda f(x)) \dd x = \cf_{\Gamma_2}(\lambda) \E(\lambda f \circ x_0)i_{x_0}^*(\varphi) \prod_{i = 1}^n \left(\int_{\B(0, \alpha)} \E(\lambda a_iy_i^2) \dd y_i\right) \]
    which concludes the proof with $A \coloneq \min\{-A_1 - N, A_2\}$.
\end{proof}

Now, the proof of Theorem~\ref{thm:FPS-mot} is complete up to showing the following lemma and proposition.

\begin{lemm}\label{lemma:Gauss-mot}
    Let $\alpha$ be a definable constant of $h[0, 0, 1]$. Let $a$ and $b$ be two definable constants of $h[1, 0, 0]^\times$ such that, for all fields $K$ of $\Field_k$, we have
    \[ \ord b_K - \ord a_K \geqslant \alpha_K \qquad\text{and}\qquad 2\ord b_K - \ord a_K \geqslant 1. \]
    Then
    \[ \int_{\B(0, \alpha)} \E(ax^2 + bx) \dd x = \int_{\B(0, \alpha)} \E(ax^2) \dd x. \]
\end{lemm}

\begin{proof}
    We use the same arguments as in the proof of Lemma~\ref{lemma:Gauss}. In particular, Axioms (R2) and (R3) of the exponential $\E$ (see~\cite[§3.1]{CluckersLoeser2010}), Axiom~(A3) of the motivic integral (see~\cite[Theorem 10.1.1]{CluckersLoeser2010}) and the change of variables theorem (see~\cite[Theorem 4.2.1]{CluckersLoeser2010}) are used.
\end{proof}

\begin{prop}\label{prop:theta-SB}
    The function $\theta$ defined by equation~\eqref{eq:theta} in the proof of Theorem~\ref{thm:FPS-mot} is a Schwartz--Bruhat function.
\end{prop}

\begin{proof}
    We use all notations of the proof of Theorem~\ref{thm:FPS-mot}. By definition of the morphism $\iota_!$, the support of the function $\theta$ is bounded. We have to show that it is locally constant. By Proposition~\ref{prop:prod-SB}, it is enough to show this property for the functions~$\iota_!((j \circ u)^*(\varphi))$ and~$\iota_!(\LL^{\ord\Jac u^{-1} \circ u})$.

    Let us work on the function $\iota_!((j \circ u)^*(\varphi))$. We set
    \[ X \coloneq \{(x, z) \in h[2n, 0, 0] \mid \ord(x - z) \geqslant \alpha^+(\varphi)\}. \]
    Let $(x, z, K)$ be a point of $X$. With Proposition~\ref{prop:SB-eval}, it is enough to show the equality
    \begin{equation}\label{eq:first-function-eval}
        \iota_!((j \circ u)^*(\varphi))(x, K)_{\S(x, z, K)} = \iota_!((j \circ u)^*(\varphi))(z, K)_{\S(x, z, K)}.
    \end{equation}
    Recall that we assumed $\alpha \geqslant \alpha^+(\varphi)$. If~$(x, K) \notin \abs{\B(0, \alpha)^n}$, then $(z, K) \notin \abs{\B(0, \alpha)^n}$ and, by Lemma~\ref{lemma:ext-zero}, the two members of the equality~\eqref{eq:first-function-eval} are equal to zero and we deduce that this equality is satisfied. We assume $(x, K) \in \abs{\B(0, \alpha)^n}$. Then~$(z, K) \in \abs{\B(0, \alpha)^n}$. Thus, by Lemma~\ref{lemma:ext-zero} and~\cite[Corollary 3.6.1]{CluckersHalupczok}, we have
    \begin{align*}
        \iota_!((j \circ u)^*(\varphi))(x, K)_{\S(x, z, K)} &= (j \circ u)^*(\varphi)(x, K)_{\S(x, z, K)} \\
        &= \varphi(j \circ u(x, K))_{\S(x, z, K)} \\
        &= \varphi(x_0 + T_K^{-1}(x), K)_{\S(x, z, K)}.
    \end{align*}
    But using the isomorphism $T \colon \B(0, \alpha)^n \to \B(0, \alpha)^n$, we can write
    \begin{align*}
        \ord(x_0 + T_K^{-1}(x) - [x_0 + T_K^{-1}(z)]) &= \ord(T_K^{-1}(x) - T_K^{-1}(z)) \\
        &\geqslant \min\{\ord T_K^{-1}(x), T_K^{-1}(z)\} \geqslant \alpha
    \end{align*}
    As $\alpha \geqslant \alpha^+(\varphi)$, if we denote by $i \colon \Omega \to h[n, 0, 0]$ the canonical inclusion morphism, the function $i_!\varphi$ is locally constant for the definable constant $\alpha$ and we finally have the equality~\eqref{eq:first-function-eval} implied by the equality
    \[ \varphi(x_0 + T_K^{-1}(x), K)_{\S(x, z, K)} = \varphi(x_0 + T_K^{-1}(z), K)_{\S(x, z, K)}. \]

    Let us work on the function $\iota_!(\LL^{\ord\Jac u^{-1} \circ u})$. By the chain rule formula for definable morphisms~\cite[Proposition 8.4.1]{CluckersLoeser2008}, we have
    \[ 0 = \ord\Jac(u^{-1} \circ u) = \ord\Jac u + (\ord\Jac u^{-1}) \circ u. \]
    But the functions $\ord\Jac u$ and $\ord\Jac u^{-1}$ are nonnegative since both morphism~$u$ and~$u^{-1}$ are analytic. We deduce that
    \[ (\ord\Jac u^{-1}) \circ u = 0 \]
    and so the function $\iota_!(\LL^{\ord\Jac u^{-1} \circ u}) = \cf_{\B(0, \alpha)^n}$ is locally constant.
\end{proof}

\subsection{A global version}
\label{ssec:global}

We can improve Theorem~\ref{thm:FPS-mot} by considering critical points as a variable. We will take the same notation as in §\ref{sec:formula-mot}. Let $\Omega$ be a definable subassignment of $\O^n$. Let $\varphi$ be a Schwartz--Bruhat function on $\Omega$ and~$f \colon \Omega \to \O$ be a globally analytic definable morphism. Let~$\Lambda$ be an unbounded definable subassignment of $h[1, 0, 0]$. We assume the following points:
\begin{enumerate}[label=(\roman*)]
    \item\label{it:points-critique} for any field $K$ of $\Field_k$, the function $f_K$ admits finitely many critical points, these critical points are nondegenerate modulo~$t$ and they belong to the polydisc~$\B(0, \alpha^-(\varphi))^n(K)$;

    \item there exists a definable morphism $R \colon \Omega[n, 0, 0] \to h[1, 0, 0]$ such that, for all fields $K$ of $\Field_k$, we have
    \[ f_K(x + y) = f_K(x) + \inner{\grad f_K(x)}{y} + \inner{R_K(x, y)y}{y} \]
    for all elements $x$ of $\Omega(K)$ and $y$ of $K\ls t^n$ with $x + y \in \Omega(K)$;

    \item\label{it:reste} there exists a definable constant $M$ of $h[0, 0, 1]$ such that, for all fields $K$ of~$\Field_k$, we have
    \[ \ord R_K(x, y) \geqslant M_K \]
    for all elements $x = (x_1, \dots, x_n)$ of $\Omega(K)$ and $y = (y_1, \dots, y_n)$ of $K\ls t^n$ with~$\ord x_i \geqslant \alpha^-(\varphi)_K$ and $\ord y_i \geqslant \alpha^+(\varphi)_K$ for all indexes $i$ in $\{1, \dots, n\}$.
\end{enumerate}

\begin{enonce}[remark]{Notation}\label{par:global-notation}
    We consider the definable subassignment
    \[ \Crit f \coloneq \{x \in \Omega \mid \grad f(x) = 0\} \subset \B(0, \alpha^-(\varphi))^n. \]
    By point~\ref{it:points-critique}, for all fields $K$ of~$\Field_k$, the set $(\Crit f)(K)$, that is
    \[ \Crit f_K \coloneq \{x \in \Omega(K) \mid \grad f_K(x) = 0\}, \]
    is finite. Consequently, we can find a definable constant $\alpha$ of $h[0, 0, 1]$ such that, for all points~$(x_0, y_0, K)$ of $(\Crit f)^2$ with $x_0 \ne y_0$, we have
    \[ (x_{0K} + t^{\alpha_K}K\ls t^n) \cap (y_{0K} + t^{\alpha_K}K\ls t^n) = \emptyset \]
    and
    \[ x_{0K} + t^{\alpha_K}K\ls t^n \subset \Omega(K). \]
    For example, we can take
    \[ \alpha_K = \max(1 + \max_{x, y \in \Crit f_K} \ord(x - y), \inf\{\beta \in \ZZ \mid x_{0K} + t^\beta K\ls t^n \subset \Omega(K)\}) \]
    for each field $K$ of $\Field_k$. We consider the definable morphism
    \[ \fonc{g}{\Lambda \times \Crit f \times \Omega}{h[1, 0, 0],}{(\lambda, x_0, x)}{\lambda f(x)} \]
    and the Presburger function given by the definable morphism
    \[ \fonc{\psi}{\Lambda \times \Crit f \times \Omega}{h[0, 0, 1],}{(\lambda, x_0, x)}{\cf_{\B(x_0, \alpha)}(x)}. \]
    Finally, let $\pi \colon \Lambda \times \Crit f \times \Omega \to \Lambda \times \Crit f$ and $\pi' \colon \Lambda \times \Crit f \times \Omega \to \Omega$ be the canonical projections.
\end{enonce}

\begin{defi}
    We consider the function
    \[ I_{f, \varphi, \alpha} \coloneq \pi_{!\Lambda \times \Crit f}(\psi\pi'^*(\varphi)\E(g)) \in \Cons(\Lambda \times \Crit f)\expe. \]
\end{defi}

Let $i \colon \Crit f \to \Omega$ be the canonical inclusion and $p \colon \Lambda \times \Crit f \to \Crit f$ be the canonical projection.

\begin{theo}
    With assumptions~\ref{it:points-critique}--\ref{it:reste} of §\ref{ssec:global} and notation~\ref{par:global-notation}, there exists a definable morphism $A \colon \Crit f \to h[0, 0, 1]$ and~$a_1, \dots, a_n \colon \Crit f \to \O^\times$ such that the equality
    \begin{multline}\label{eq:fps-param}
        I_{f, \varphi, \alpha}\cf_\Gamma = \E((\lambda, x_0) \longmapsto \lambda f(x_0))(p \circ i)^*(\varphi) \\
        \prod_{i = 1}^n \mu_{\Lambda \times \Crit f}(\cf_{\Lambda \times \Crit f \times \B(0, \alpha)}\E((\lambda, x_0, x) \longmapsto \lambda a_i(x_0)x^2))\cf_\Gamma
    \end{multline}
    holds in the ring $\Cons(\Lambda \times \Crit f)\expe$ where
    \[ \Gamma \coloneq \{(\lambda, x_0) \in \Lambda \times \Crit f \mid \ord\lambda \leqslant A(x_0)\}. \]
\end{theo}

\begin{proof}
    We follow the same argument as in the proof of Theorem~\ref{thm:FPS-mot} working relatively to the definable subassignment $\Crit f$. On can check in a straightforward way that the quantities $A$, $\beta$, $a_1$, \dots, $a_n$ depend definably of critical points.
\end{proof}

\begin{rema}
    Using evaluation as §\ref{sssec:eval}, we can interpret formula~\eqref{eq:fps-param}. Indeed, for any point $(\lambda, x_0)$ of $\Lambda \times \Crit f$ with $\ord\lambda \leqslant A(x_0)$, we get
    \[ \int_{\B(x_0, \alpha)} \varphi(x)\E(\lambda f(x)) \dd x = \E(\lambda f(x_0))\varphi(x_0) \prod_{i = 1}^n \int_{\B(0, \beta(x_0))} \E(\lambda a_i(x_0)x^2) \dd x \]
    in the ring $\Cons(\{(\lambda, x_0)\})\expe$.

    The main advantage of formula~\eqref{eq:fps-param} is that it does not involve morphisms of the form $i_{x_0}^*$ and, in some way, evaluations of constructible functions. To be more precise, comparing with Theorem~\ref{thm:FPS-mot}, we have replaced these evaluations by the definable morphism $p \circ i \colon \Lambda \times \Crit f \to \Omega$.
\end{rema}

\section*{Acknowledgements}

The author would like to thank his advisor Michel~Raibaut for its guidance, Arthur~Forey for interesting discussions and the referee for the meticulous review of the paper and useful remarks.

\printbibliography

\end{document}